\documentclass[12pt]{article}
\usepackage{amsthm,amsfonts,amsmath,amssymb}

\newcommand{\diff}{{\rm d}}
\newcommand{\eq}{\begin{equation}}
\newcommand{\en}{\end{equation}}

\newcommand{\re}[1]{\mbox{(\ref{#1})}}

\newcommand{\prob}{\mathbb P}

\newcommand{\ex}{\mathbb E}
\newcommand{\nints}{\mathbb N}
\newcommand{\Nat}{\Bbb N}

\def\endpf{\hfill $\Box$ \vskip0.5cm}
\def \proof{\noindent{\it Proof.\ }}

\newcommand{\giv}{\,|\,}
\newtheorem{theorem}{\large Theorem}

\newtheorem{corollary}[theorem]{\large  Corollary}
\newtheorem{lemma}[theorem]{\large  Lemma}

\begin{document}

\title{Poisson representation of a Ewens fragmentation process
\thanks{Research supported in part by N.S.F. Grant DMS-0405779}
}
\author{Alexander Gnedin\thanks{Mathematical Institute, Utrecht University, The Netherlands; e-mail gnedin@math.uu.nl}
\hspace{.2cm}
and 
\hspace{.2cm}
Jim Pitman\thanks{Department of Statistics, University of California, Berkeley, USA; e-mail pitman@stat.Berkeley.EDU} 
\\
\\
\\
\\
}
\date{
\today
\\
}
\maketitle

\maketitle
\begin{abstract}
\noindent
A simple explicit construction is provided of a partition-valued
fragmentation process whose distribution on partitions of $[n]=\{1,\ldots,n\}$ at time $\theta \ge 0$ 
is governed by the Ewens sampling formula with parameter $\theta$.
These partition-valued processes are exchangeable and consistent, as $n$ varies.
They can be derived by uniform sampling from a corresponding mass fragmentation process 
defined by cutting a unit interval at the points of a Poisson process with intensity
$\theta x^{-1} \diff x$ on ${\mathbb R}_+$, arranged to be intensifying as $\theta$ increases.
\end{abstract}

\vskip0.5cm

\noindent

\section{Introduction}

There has been much recent interest in models for random processes of fragmentation
and coagulation: see Chapter 5 of \cite{csp} and the recent book
\cite{bertoin-coag-frag}. 
Mekjian and others \cite{chasemek94,lm91,mekj91,ml91} have considered
Ewens partitions with parameter $\theta$ as a model for fragmentation
phenomena, with the intuitive notion that increasing $\theta$ corresponds
to further fragmentation. But it does not seem obvious how to
construct a nice Markovian fragmentation process corresponding to this idea.

It was pointed out in \cite{csp} that it is possible to construct a sequence of 
partition-valued processes $(\Pi_{n,\theta}\,,\, \theta \ge 0)$ ($n = 1,2, \ldots$)
with the following properties:
\begin{itemize}
\item (Ewens distribution)
$\Pi_{n,\theta}$ is for each $n = 1,2, \ldots$ and $\theta \ge 0$
a random partition of the set $[n]:= \{1,2, \ldots, n \}$,
with distribution determined by the following formula:
for each composition $(n_1, \ldots, n_k)$  of $n = \sum_{i=1}^k n_i\,$,
and each partition $\pi$ of $[n]$ into $k$ blocks of sizes $n_1, \ldots, n_k$,
\eq
\label{eppf}
\prob (\Pi_{n,\theta} = \pi ) = \frac{\theta^{k-1}} { [\theta + 1 ]_{n-1} } \prod_{i = 1}^k (n_i - 1)!
\en
where 
$[x]_m := x(x+1)\cdots(x+m-1)$ is the Pochhammer factorial;
\item (fragmentation) $\Pi_{n,\theta}$ is a refinement of $\Pi_{n,\phi}$ if $\theta > \phi$, for all $n \ge 1$,
that is each block of $\Pi_{n,\phi}$ is some union of blocks of  $\Pi_{n,\theta}$; 
\item (consistency) for each $n < m$ 
a process with the same distribution as $(\Pi_{n,\theta}\,,\, \theta \ge 0 )$ is
obtained by restriction of $(\Pi_{m,\theta}\,,\, \theta \ge 0 )$ to $[n]$;

\item (exchangeability) for each $n$ the law of $(\Pi_{n,\theta}\,,\, \theta \ge 0)$
is invariant under permutations of the set $[n]$.

\end{itemize}
We call a sequence of partition-valued processes $(\Pi_{n,\theta}\,,\,\theta\geq 0)$ ($n=1,2,\ldots$) 
with these properties a  
{\em family of Ewens fragmentations}.
One family of Ewens fragmentations associated with Kingman's coalescent was analysed in \cite{beresp06}.

\par Some general theory \cite{jpse96cmc} implies that any such family of processes 
can be defined in the strong sense consistently (i.e. so that $\Pi_{m,\theta}|_{[n]}=\Pi_{n,\theta}$ for $m>n$)
on a single probability space by means of uniform sampling 
of points engaged in a process of 
fragmentation of a total mass  $1$ into a countable collection of submasses with sum $1$, with more and more
refined splitting of the submasses as the time parameter $\theta$ increases.
See \cite[Chapter 3]{bertoin-coag-frag} and \cite[Chapter 5]{csp} for further background and references.

\par The problem was posed in \cite{csp} of characterising the dynamics of a  family of Ewens 
fragmentations, preferably in a Markovian way. 
For applications, it is desirable to have a model which can easily 
be simulated for modest values of $n$.  But previous efforts fall short in this respect.
In this note we partly solve this problem by constructing 
a new family of Ewens fragmentations. Our family is not  Markovian,
but it enjoys the Markov property 
and follows a very simple transition rule
when viewed as a fragmentation process in the extended space of {\it ordered}
partitions.
This simplification by passing to an ordered structure extends our previous work on 
regenerative partitions and their relatives  \cite{RPS, RCS,  selfsim}.

\section{Construction}

Recall that a {\em composition} of $n$ is a sequence of positive integers
$(n_1, \ldots, n_k)$ with sum $n$. We regard a composition of $n$ as a way
of distributing $n$ unlabelled balls in an ordered sequence of $k$ non-empty
boxes, with $n_i$ balls in the $i$th box.  
A composition of $n$ is also conveniently encoded by the binary sequence of $0$'s and $1$'s obtained
by concatenating subsequences of the form $1$, $10$, $100$, $\ldots$, where the
$i$th subsequence in the concatenation has length $n_i$. So the symbols $1$ occur at places
$1, n_1+1, n_1 +n_2 +1, \ldots \sum_{i = 1}^{k-1} n_i + 1$.
Using a particular composition $(3,4,1)$ of $8$ for illustration,  
the balls-in-boxes picture is suggested by the notation $[000] \, [0000] \, [0]$.
The binary representation is obtained by replacing each $[0$ by a $1$ and deleting each $]$ to obtain
the sequence $10010001$.
Let $x$ and $y$ be two compositions of $n$, each represented as a binary sequence, say
$x = (x_i, 1 \le i \le n )$ and $y = (y_i, 1 \le i \le n )$.
Say $x$ is a {\em refinement} of $y$ if $x_i \ge y_i$ for every $i$. 
In terms  of the balls-in-boxes picture, $x$ is derived from $y$ by splitting boxes
into sub-boxes, and in terms of the binary representation
$x$ is derived by switching some $0$'s to $1$'s.
For instance, $11010101$ is a refinement of  $10010001$.

Given a stochastic process $(C_{n,\theta}\,,\, \theta \ge 0 )$, with values in compositions of
$n$ we define an {\em associated partition-valued process} $(\Pi_{n,\theta}, \theta \ge 0 )$ by first 
assigning each of the $n$ places for a ball in the balls-in-boxes representation a number in $[n]$ according to 
a uniform random permutation of $[n]$ independent of $(C_{n,\theta}\,,\, \theta \ge 0 )$, 
thus obtaining an {\it ordered partition} of $[n]$,
then ignoring the
order of the boxes to obtain a partition of $[n]$.
If $(C_{n,\theta}\,,\, \theta \ge 0 )$ is refining as $\theta$ increases then
the associated partition-valued process $(\Pi_{n,\theta}\,,\, \theta \ge 0 )$ is a fragmentation process
whose law is invariant under permutations of $[n]$.
The process $(\Pi_{n,\theta}\,,\, \theta \ge 0 )$ then describes a process of randomly splitting up a
collection of balls labelled by $[n]$ into an unordered collection of boxes.

\begin{theorem}\label{1}
Let $\Theta_j$ for $j = 1,2, \ldots$ be a sequence of independent random variables with distributions
\eq
\label{thj}
\prob (\Theta _j \le \theta ) = \theta/(\theta + j - 1) , ~~~\theta \ge 0
\en
(where $0/0=1$).
Let $C_{n,\theta}$ for $n = 1,2, \ldots$ be the random composition of $n$ whose binary representation
is the sequence of indicator variables $1(\Theta_j \le \theta)$ for $1 \le j \le n$.
Then the sequence $(\Pi_{n,\theta}\,,\, \theta \ge 0)$ of partition-valued process 
associated with $(C_{n,\theta}\,,\, \theta \ge 0)$ defines a
family of Ewens fragmentations.
\end{theorem}

That $\Pi_{n,\theta}$ has the Ewens distribution \re{eppf} can be read from the known result 
\cite{ABT, selfsim, csp} 
that in the binary representation of the composition of $n$ derived from the block sizes of a Ewens partition of $[n]$ in 
reversed size-biased order, 
the digits are independent Bernoulli variables with parameters
$\theta/(\theta + j - 1)$ as in \re{thj}.
The device \re{thj} with independent variables $\Theta_j$ is then the simplest way to make these indicators
simultaneously for all $j$ and $n$ to be increasing in $\theta$, which is all that is needed to make
$(\Pi_{n,\theta}\,,\, \theta \ge 0)$ a Ewens fragmentation.
What is much less obvious is the consistency of these  processes for various $n$. To put this another way,
if in the process of splitting a set of $m$ balls according to the indicators $1(\Theta_j \le \theta)$ for $1 \le j \le m$
we pass to the balls-in-boxes picture and just observe the splitting process restricted to a uniformly chosen
random subset of $n < m$ balls, this sub-process is identical in distribution to the process of 
splitting of the first $n$ balls using the indicators $1(\Theta_j \le \theta)$ for $1 \le j \le n$.
This sampling consistency property of compositions, which is so intuitive in the balls-in-boxes picture, is
quite painful to express entirely in the binary encoding. We circumvent that difficulty  by deriving 
 consistency from the Poisson representation of the corresponding mass fragmentation model, 
which is introduced in the next section.
See also \cite{selfsim} for a more extensive discussion of such consistency properties of 
partition structures derived by random sampling from self-similar random sets,
like the self-similar Poisson process in the next section
or the zero set of a Brownian motion.

\section{Poisson representation of the mass fragmentation}

Let $(T_i,V_i)$ be a listing of the points of a homogeneous Poisson point process on the positive
quadrant ${\mathbb R}_+^2$ with rate 1 per unit area.
For each fixed $\theta >0$, the random countable set
$$
Z_\theta := \{ T_i : 0 < T_i < 1 \mbox{ and } V_i \le \theta/T_i \}
$$
is then the set of points of a Poisson process on $[0,1]$ with intensity $\theta t ^{-1} \diff t $.
To orient the reader, we start by recalling some well known properties of $Z_\theta$ and
the induced random composition \cite{ABT,  selfsim, tsf, csp}.
\begin{enumerate}
\item[(i)]
Let $Y_{1,\theta} > Y_{2,\theta} > \cdots$ be the points of $Z_\theta$ in decreasing order.
Then
$$Y_{j,\theta} = \prod_{i = 1}^j W_{i,\theta}$$ where the $W_{i,\theta}$ ($i=1,2,\ldots$) are independent
and identically distributed random variables with beta$(\theta,1)$ distribution. 
\item[(ii)]
If $P_{j,\theta}$ is the length of the $j$th component interval of $[0,1] \backslash Z_\theta$,
working from right to left, that is $P_{j,\theta} = Y_{j-1,\theta }  - Y_{j,\theta}$ where $Y_{0,\theta}:= 1$, then
$$P_{j,\theta} = (1 - W_{j,\theta}) \prod_{1 = 1}^{j-1} W_{i,\theta}$$ for $W_{i,\theta}$ 
i.i.d.  beta$(\theta,1)$ as before.
The distribution of this random discrete probability distribution $(P_{j,\theta}\,,\, j \ge 1 )$ is known as the
GEM$(\theta)$ distribution.
\item[(iii)]
The distribution of the decreasing rearrangement of $(P_{j,\theta}\,,\, j \ge 1 )$ is the 
{\em Poisson-Dirichlet distribution}
with parameter $\theta$.
\item[(iv)]
Let $U_1, U_2, \ldots, $ be a sequence of independent uniform $[0,1]$ variables, independent of $Z_\theta$, and
define a random partition $\Pi_{\infty,\theta}$ of the set $\nints$ of
positive integers to be the collection of equivalence classes
for the random equivalence relation:\, $i \sim^\theta j$ if and only if either $i=j$ or both $U_i$ and $U_j$ fall in the
same component interval of $[0,1] \backslash Z_\theta$. 
Then $\Pi_{\infty,\theta}$ is an exchangeable random partition of the infinite set $\nints$,
whose restriction $\Pi_{n,\theta}$ to $[n]$ is a Ewens partition governed by 
(\ref{eppf}) for each $n$.
\item[(v)]
Let $U_{n,j}$ be the $j$th smallest value among $U_1, \ldots, U_n$, and let $X_{n,j}(\theta)$ be the indicator of the
event that $U_{n,j}$ is the least value among those of the $n$ values 
which fall in some component interval of $[0,1]\backslash Z_\theta$.
Then for each fixed $n$ and $\theta$, the $X_{n,j}(\theta)$ for $j = 1,2, \ldots , n$ are independent, with
\eq
\label{xj}
\prob ( X_{n,j}(\theta) = 1 ) =   \theta/(\theta + j -1 ).
\en
More precisely, the sequence $(X_{n,j}(\theta), 1 \le j \le n)$ is the binary encoding of
a random composition $C_{n,\theta}$ of $n$ which is a particular ordering of the sizes of blocks of $\Pi_{n,\theta}$.
If the blocks of $\Pi_{n,\theta}$ are, say, $B_1,\ldots,B_k$, then the
sizes of these blocks appear in the composition
in increasing order of values of $\min_{i \in B_k} U_i\,$\footnote{Strictly speaking, this 
defines $C_{n,\theta}$  in terms of $U_j$'s and $Z_\theta$, rather than through $\Pi_{n,\theta}$.
Conditionally given $\Pi_{n,\theta}=\{B_1,\ldots,B_k\}$ the composition has the same distribution
as the sequence $\#B_j$ arranged by decrease of minimal elements $\min B_j$.
Conditionally given the induced partition 
$\{\#B_j\,,\,j\leq k\}$ of integer $n$,
this arrangement is  the inverse size-biased ordering of the block sizes.}.
\end{enumerate}

An immediate consequence of the above construction of $Z_\theta$ from a Poisson process in the positive quadrant is that the random 
set $Z_\theta$ increases as $\theta$ increases. Consequently, the various quantities 
introduced above describe a process of fragmentation of 
 $[0,1]$ into subintervals. 
In particular, the partition valued process $(\Pi_{n,\theta}\,,\, \theta \ge 0)$ is refining as $\theta$ increases,
and each of the processes $(X_{n,j}(\theta)\,,\, \theta \ge 0 )$ is increasing as $\theta$ increases.
Consistency and exchangeability of the processes $(\Pi_{n,\theta}\,,\, \theta \ge 0)$ 
 are obvious from (iv).

\par The proof of the consistency property claimed in Theorem \ref{1} is completed by the following lemma,
which shows that in this Poisson setup the indicator variables $(X_{n,j}(\theta), 1 \le j \le n)$ in the 
binary expansion of the composition $C_{n,\theta}$  
associated with 
the natural ordering of blocks of 
$\Pi_{n,\theta}$ (as in (v)) can be derived from independent variables $(\Theta_{n,j}\,,\, 1 \le j \le n)$ 
with the same distribution as $(\Theta_{j}\,,\, 1 \le j \le n)$ in Theorem \ref{1}.

\begin{lemma}
Let $\Theta_{n,j} : = \inf \{ \theta : X_{n,j} (\theta ) = 1 \}$, so that $X_{n,j}(\theta)$ may be represented as
$$
X_{n,j}(\theta) = 1 (\Theta_{n,j } \le \theta ).
$$
Then for each fixed $n$ the $\Theta_{n,j}$ are independent random variables with
$$
\prob ( X_{n,j}(\theta)  = 1 ) = \prob  (\Theta_{n,j } \le \theta ) = \theta/(\theta + j - 1)\,,~~~~(j=1,\ldots,n).
$$
\end{lemma}
\proof
Fix $\theta_j > 0$ for $j = 2, \ldots , n $.
Observe that the event $(\Theta_{n,j} > \theta_j)$ occurs if and only if there is no point $(T_i,V_i)$ of the
Poisson process with  $T_i \in [U_{n,j-1}, U_{n,j}]$ and $V_i \le \theta/T_i$ .
Therefore
$$
\prob \left( \cap_{j = 2 }^n (\Theta_{n,j} > \theta_j ) \,|\,U_{n,j}, 1 \le j \le n \right) 
= \exp \left( - \sum_{j = 2}^n \int_{U_{n,j-1}} ^{U_{n,j}} \frac {\theta_j}{t}\, \diff t \right)
= \prod_{j = 2}^n \left( \frac{U_{n,j-1} } { U_{n,j } } \right)^{\theta_j}
$$
and hence
$$
\prob ( \cap_{j = 2 }^n (\Theta_{n,j} > \theta_j ) 
= \ex \left [ \prod_{j = 2}^n \left( \frac{U_{n,j-1} } { U_{n,j } }\right)^{\theta_j} \right]
= \prod_{j = 2}^n  \frac {j - 1 }{ j - 1 + \theta_j}\,\,,
$$
because the ratios ${U_{n,j-1} }/ { U_{n,j }} $ are independent with beta$(j-1,1)$ distributions, $2 \le j \le n$.
\endpf

\par As before, let $(C_{n,\theta}\,,\,\theta\geq 0)$ 
be the process of refining compositions of $n$, defined either through indicators as in  Theorem \ref{1},
or by means of the Poisson construction as in (v) above.
Immediately from the definition, we have:

\begin{corollary}\label{3}
$(C_{n,\theta}, \theta \ge 0)$ is a Markov process whose
inhomogeneous transition rates  are determined by the rule:
if at time $\theta$ the state is the composition of $n$ encoded by
some binary sequence starting  with   $1$,
each $0$ is switching to $1$ at rate $1/(\theta + j-1)$, 
where $j$ is the place of this $0$ in the sequence, while all other transition 
rates are trivial.
\end{corollary}
\proof Indeed,
$$\prob(\Theta_j\in [\theta\,,\,\theta+\diff\theta]\giv\Theta_j>\theta)={\left({\theta\over\theta+j-1}\right)' \diff\theta
\over 1-{\theta\over \theta+j-1}}={\diff\theta\over \theta+j-1}\,.$$
\endpf

\par If we view  $\Pi_{\infty,\theta}$ together with the total ordering of  the blocks, as induced by
the natural order of intervals (recall (iv)), we obtain an exchangeable ordered partition $\Pi_{\infty,\theta}^*$  
of $\Nat$.
Let $\Pi_{n,\theta}^*$ be its restriction to $[n]$. The law of $\Pi_{n,\theta}^*$ is given by  
the ordered version of the Ewens sampling formula (\ref{eppf}): 
$$\prob(\Pi_{n,\theta}^*=\pi^*)={\theta^{k-1}\over[\theta+1]_{n-1}}\prod_{j=1}^k {n_j!\over n_1+\cdots+n_j}\,$$
for every $\pi^*$ ordered partition of $[n]$ with block sizes $(n_1,\ldots,n_k)$.
The process $(\Pi_{n,\theta}\,,\,\theta\geq 0)$ is Markovian for every $n$.
The transition mechanism of $(\Pi_{n,\theta}^*\,,\,\theta\geq 0)$ is determined by that of 
$(C_{n,\theta}\,,\,\theta\geq 0)$ and the following allocation rule\footnote{Which is common for all
exchangeable fragmentation processes.}: each time a block $B_j$ of size $a$ splits
in two fragments of sizes $\xi$ and $\eta$, all ${a\choose \xi}$ possible allocations of the elements of 
$B_j$ among the offspring fragments are equally likely.

\section{Further properties}

In principle, the finite dimensional distributions of 
$(\Pi_{n,\theta}\,,\, \theta \ge 0)$ may be determined by some summations  of
probabilities determined by the Markov process $(C_{n,\theta}, \theta \ge 0)$. But such formulas appear to be of
limited value.  It appears that the process $(\Pi_{n,\theta}, \theta \ge 0)$ is not Markovian.

\vskip0.5cm
\noindent
{\bf Proof that for $n>2$ the fragmentation process is not Markovian} 
Consider the random time $\widehat{\Theta}_n$ of the first split, 
that is
$$
\widehat{\Theta}_n 
 = \min_{2 \le j \le n} \Theta_{n,j} 
= \inf \{ \theta : \Pi_{n,\theta } \neq \{[n]\}\},
$$
where $\Pi_{n,0}=\{[n]\}$ is the initial partition with a single  block.
To show that the Markov property of the fragmentation process
$(\Pi_{n,\theta}, \theta \ge 0)$
does not hold for every $n>2$ we shall focus on the conditional probability
$$Q(t):
=\prob(\Pi_{n,\theta}=\lambda\giv\Pi_{n,\phi}=\lambda\,,\,\Pi_{n,t}=\lambda)
=\prob(\Pi_{n,\theta}=\lambda\giv\Pi_{n,\phi}=\lambda\,,\,\widehat{\Theta}_n<t),$$
where $\phi,\theta$ are considered as parameters, $0<t<\phi<\theta$, and 
$\lambda$ is the partition of $[n]$ in two blocks $\{1\}$ and $\{2,\ldots,n\}$.
To disprove the Markov property it is sufficient to show that $Q(t)$ is not constant as $t$ varies.

\par Note that $\Pi_{n,\phi}=\lambda$ is only possible when the composition $C_{n,\phi}$ assumes either
the value $1100\ldots0$ or the value $100\ldots01$, and conditionally given either of these values 
$\Pi_{n,\phi}$ equals $\lambda$ with probability $1/n$
(as a consequence of exchangeability).
Working out this dichotomy,

\begin{align*}
\prob(\Pi_{n,\phi}=\lambda,\widehat{\Theta}_n<t)~=~~~~~~~~~~~~~~~~~~~~~~~~~~~~~~~~~~~~~~~~~~~~~~~~~~~~~~~~~~~~~~~\\
\prob(\Theta_{n,2}<t,\Theta_{n,n}>\phi\,;\,\Theta_{n,3}>\phi,\ldots, \Theta_{n,n-1}>\phi)\,{1\over n}~+~~~~~~~~~~
~~~~~~~~~~~~~~\\
\prob(\Theta_{n,2}>\phi,\Theta_{n,n}<t\,;\,\Theta_{n,3}>\phi,\ldots, \Theta_{n,n-1}>\phi)\,{1\over n}~=~~~~~~~~~~~~~~~
~~~~~~~~~\\
{t\over t+2-1}\left(1-{\phi\over \phi+n-1}\right)P(\phi)\,{1\over n}~+~
\left( 1- {\phi\over\phi+2-1}\right) {t\over t+n-1}\,P(\phi)\,{1\over n}~=\\
\left({t(n-1)\over (t+1)(\phi+n-1)}+{t\over (\phi+1)(t+n-1)}\right)P(\phi){1\over n}\,,~~~~~~~~~~~~~~~~~~~~~~~~~~~~
\end{align*}
where for shorthand $P(\phi):=\prob(\Theta_{n,3}>\phi,\ldots , \Theta_{n,n-1}>\phi)$.
Noting the inclusion 
 $$ \{\Pi_{n,\theta}=\lambda,\widehat{\Theta}_n<t\}\subset\{\Pi_{n,\phi}=\lambda,\widehat{\Theta}_n<t\}$$
and applying the above formula to the event on the left-hand side, we compute
$$Q(t)=
{\prob(\Pi_{n,\theta}=\lambda\,,\,\widehat{\Theta}_n<t)\over
\prob(\Pi_{n,\phi}=\lambda\,,\,\widehat{\Theta}_n<t)}=
\left({t(n\theta+2n-2)+(n-1)^2(\theta+1)+\theta+n-1\over t(n\phi+2n-2)+(n-1)^2(\phi+1)+\phi+n-1}\right)
{P(\phi)\over P(\theta)}\,.$$
This does not depend on $t$  if and only if
$${n\theta+2n-2\over n\phi+2n-2}={(n-1)^2(\theta+1)+\theta+n-1\over (n-1)^2(\phi+1)+\phi+n-1}\,,$$
or, equivalently, if and only if the polynomial
\begin{align*}
(n\theta+2n-2)( (n-1)^2(\phi+1)+\phi+n-1)~=~~~~~~~~~~~~~~~~~~~~\\
n(n^2-2n+2)\theta\phi+n^2(n-1)\theta+2(n-1)(n^2-2n+2)\phi
+2n(n-1)^2
\end{align*}
is symmetric in $\phi$ and $\theta$. To maintain symmetry we must have
$$n^2(n-1)=2(n-1)(n^2-2n+2),$$
which forces positive $n$ to be either $1$ or $2$. 
Thus for $n>2$ the partition-valued process is not Markovian (while it is trivially Markovian for
$n=1$ or $2$).

\vskip0.5cm
\noindent
{\bf Transition rates of the Ewens fragmentation}
Given the value $\pi$ of $\Pi_{n,\theta}\,$, the composition $C_{n,\theta}$ can be recovered by
arranging the sequence of block sizes of $\pi$ in the reversed size-biased order.
This property taken together with Corollary \ref{3} allows to compute the transition rates.
To illustrate the method, suppose that at time $\theta$ the partition
$\Pi_{n,\theta}$ is in  state $\pi$ with  block sizes 
$\{a,b\}$, and let $\sigma$ be some nontrivial refinement of $\pi$ with  
block sizes 
$\{\xi,\eta,b\}$, so $\xi+\eta=a$.
Suppose first that $a\neq b$, $\xi\neq\eta$.
Then $C_{n,\theta}=(a,b)$ with probability $b/(a+b)$, and $C_{n,\theta}=(b,a)$ with probability 
$a/(a+b)$. Inspecting all possibilities we see that the fragmentation process jumps from 
$\pi$ to $\sigma$ at rate
$${b\over a+b}\left({1\over \theta+\xi}+{1\over \theta+\eta}\right){a\choose\xi}^{-1}+
{a\over a+b}\left({1\over \theta+\xi+b}+{1\over \theta+\eta+b}\right){a\choose\xi}^{-1},$$
where the binomial coefficients accounts for the number of ways 
of allocating the elements of the splitting block among two new fragments.
A minute thought shows that the formula is still valid in the case $a=b$;
but if $\xi=\eta$ the above expression should be halved. 

\par In principle, there exists a Markovian family of Ewens fragmentations,
with the same 
transition rates as that of $(\Pi_{n,\theta})$'s. However, this seems to be of little use,
because the 
formulas for these rates become increasingly complicated when the number of blocks 
grows.

\vskip0.5cm
\noindent
{\bf Comparison with the Ewens fragmentation derived from Kingman's coalescent}
Another family of Ewens fragmentations was derived in \cite{beresp06} from
Kingman's coalescent tree. These fragmentations are not Markovian 
in the proper sense, and no extended Markov property for them is known.
We show next that the partition-valued process in \cite{beresp06} is different from 
the process constructed in this paper.

\par As before,
consider the  time   of the first split $\widehat{\Theta}_n$.
Let $I_n + 1$ be the almost surely unique index $j$ which makes
$\widehat{\Theta}_n  = \Theta_{n,j}$.  Then the split at time $\widehat{\Theta}_n$ creates
a partition $\Pi_{n,\widehat{\Theta}_n}$ with two blocks of sizes $I_n$ and $n-I_n$.
Conditioning on $\widehat{\Theta}_n$ gives 
\eq
\label{thn}
\prob( I_n  = i ) = \int_0^\infty 
\frac{ \prob ( \Theta_{n,i+1} \in \diff \theta) } { \prob ( \Theta_{n,i+1} > \theta)}
\prod_{2 \le j \le n } \prob ( \Theta_{n,j} > \theta)
\en
which simplifies to
\eq
\label{thn1}
\prob( I_n  = i ) = 
(n - 1)!
\int_0 ^\infty
\frac{ \diff \theta } { (\theta + i ) [ \theta + 1 ]_{n-1} }\,.
\en
In particular, for $n = 4$ this gives
\begin{eqnarray*}
\label{thn2}
\prob( I_4  = 1 ) = 6 \left( - \log 2 + \frac{1}{4} \log 3  + \frac{1}{2} \right),\\
\label{thn3}
\prob( I_4  = 2 ) = 6 \left( \frac{1}{2} \log 3  - \frac{1}{2} \right),\\
\label{thn4}
\prob( I_4  = 3 ) = 6 \left( \log 2  - \frac{3}{4} \log 3  - \frac{1}{6} \right).
\end{eqnarray*}
Compare the second of these evaluations 
with the corresponding formula
in \cite[Section 7.1]{beresp06}
to see that this Ewens fragmentation process evolves differently to the
Ewens fragmentation derived from Kingman's coalescent, which has a different distribution on
partitions of $[4]$ with two blocks at the time of the first split. 
Neither of these distributions is that of $\Pi_{4,\theta}$ given that
this partition has exactly two blocks, even though this conditional 
distribution does not depend on $\theta$.
This common conditional distribution is the Gibbs distribution on partitions of
$[4]$ into $2$ blocks which assigns probability proportional to 
$(n_1 - 1)!(n_2-1)!$ to each partition of $[4]$ into two blocks of
sizes $n_1$ and $n_2$.

\vskip0.5cm

It is still an open question for which $n$ there exists a discrete
time fragmentation process on partitions of $[n]$ whose distribution
at time $k$ is the distribution on partitions of $[n]$ into $k$ blocks
which assigns each such partition into blocks of sizes $\{n_1,\ldots,n_k\}$
a probability proportional to $\prod_{i = 1}^k (n_i - 1)!$.


\begin{thebibliography}{99}

\bibitem{ABT} Arratia, R., Barbour, A.D. and Tavar{\'e}, S. (2003) {\it Logarithmic combinatorial structures:
A probabilistic approach}, Vol. 1 of {\it EMS Monographs in Mathematics}, European Mathematical Society 
Publishing House, Z{\"u}rich. 


\bibitem{beresp06}
 Berestycki, N., and Pitman,  J. (2006)
\newblock {Gibbs distributions for random partitions generated by a
  fragmentation process}.

\bibitem{bertoin-coag-frag}
 Bertoin, J. (2006)
\newblock {\em Random fragmentation and coagulation processes}.
\newblock Cambridge University Press.


\bibitem{chasemek94}
Chase, K.C. and  Mekjian, A.Z.  (1994)
\newblock Nuclear fragmentation and its parallels.
\newblock {\em Phys. Rev. C} {\bf  49} 2164--2176.

\bibitem{jpse96cmc} 
 Evans, S.~N. and Pitman, J. (1998)
\newblock {Construction of Markovian coalescents}.
\newblock {\em Ann. Inst. Henri Poincar{\'e}} {\bf 34} 339--383.

\bibitem{tsf} Gnedin, A. (2004) Three sampling formulas,
{\em Combin.  Probab.  Comput.} {\bf 13} 185-193.


\bibitem{RPS}
Gnedin, A. and  Pitman, J. (2004)
\newblock Regenerative partition structures.
\newblock {\em Electron. J. Combin.} {\bf 11} (2), Research Paper 12, 21 pp.
  (electronic).

\bibitem{RCS}
 Gnedin, A. and  Pitman, J. (2005)
\newblock Regenerative composition structures.
\newblock {\em Ann. Probab.}, {\bf 33} 445--479.

\bibitem{selfsim}
Gnedin, A. and Pitman, J. (2005)
\newblock {Self-similar and Markov composition structures}.
\newblock In A.~A. Lodkin, editor, {\em {Representation Theory, Dynamical
  Systems, Combinatorial and Algorithmic Methods. Part 13}}, volume {\bf 326} of {\em
  Zapiski Nauchnyh Seminarov POMI} 59--84. 

\bibitem{lm91}
 Lee, S.~J. and Mekjian,  A.~Z. (1992)
\newblock Canonical studies of the cluster distribution, dynamical evolution,
  and critical temperature in nuclear multifragmentation processes.
\newblock {\em Phys. Rev. C} {\bf 45} 1284--1310.

\bibitem{mekj91}
 Mekjian, A.~Z. (1991)
\newblock Cluster distributions in physics and genetic diversity.
\newblock {\em Phys. Rev. A} {\bf 44} 8361--8374.

\bibitem{ml91}
 Mekjian, A.~Z. and  Lee, S.~J. (1991)
\newblock {Models of fragmentation and partitioning phenomena based on the
  symmetric group $S_n$ and combinatorial analysis}.
\newblock {\em Phys. Rev. A}  {\bf 44} 6294--6312.

\bibitem{csp}
Pitman, J. (2006)
\newblock {\em {Combinatorial Stochastic Processes}}.
{ Lecture
  Notes in Mathematics}, Springer.
\end{thebibliography}

\def\cprime{$'$} \def\polhk#1{\setbox0=\hbox{#1}{\ooalign{\hidewidth
  \lower1.5ex\hbox{`}\hidewidth\crcr\unhbox0}}} \def\cprime{$'$}
  \def\cprime{$'$} \def\cprime{$'$}
  \def\polhk#1{\setbox0=\hbox{#1}{\ooalign{\hidewidth
  \lower1.5ex\hbox{`}\hidewidth\crcr\unhbox0}}} \def\cprime{$'$}
  \def\cprime{$'$} \def\polhk#1{\setbox0=\hbox{#1}{\ooalign{\hidewidth
  \lower1.5ex\hbox{`}\hidewidth\crcr\unhbox0}}} \def\cprime{$'$}
  \def\cprime{$'$} \def\cydot{\leavevmode\raise.4ex\hbox{.}} \def\cprime{$'$}
  \def\cprime{$'$} \def\cprime{$'$} \def\cprime{$'$} \def\cprime{$'$}
  \def\cprime{$'$} \def\cprime{$'$} \def\cprime{$'$} \def\cprime{$'$}
  \def\cprime{$'$}

\end{document}